\documentclass[12pt]{amsart}

\usepackage{amssymb}
\usepackage[initials]{amsrefs}
\usepackage[all, poly]{xy}

\DeclareMathOperator{\charac}{char}
\DeclareMathOperator{\image}{im}
\DeclareMathOperator{\Fun}{Fun}

\DeclareMathOperator{\coev}{coev}
\DeclareMathOperator{\ev}{ev}
\DeclareMathOperator{\id}{id}
\DeclareMathOperator{\Hom}{Hom}

\DeclareMathOperator{\Tr}{Tr}

\swapnumbers

\theoremstyle{plain}
\newtheorem*{lemma}{Lemma}
\newtheorem*{proposition}{Proposition}
\newtheorem*{theorem}{Theorem}
\newtheorem*{corollary}{Corollary}

\numberwithin{equation}{subsection}

\begin{document}

\title{Quiver varieties and Lusztig's algebra}

\author{Anton Malkin}
\address{Department of Mathematics, University of Illinois,
273 Altgeld Hall, MC-382, 1409 West Green Street, Urbana, 
Illinois 61801, USA}
\email{malkin@math.uiuc.edu}

\author{Viktor Ostrik}
\address{
Department of Mathematics, University of Oregon, Eugene,
Oregon 97403-1221, USA}
\email{vostrik@darkwing.uoregon.edu}

\author{Maxim Vybornov}
\address{Department of Mathematics and Statistics,
Boston University, 111 Cummington Street, Boston,
Massachusetts 02215, USA}
\email{vybornov@bu.edu}

\begin{abstract}
We study preprojective algebras of graphs and their relationship
to module categories over representations of quantum $SL(2)$.
As an application, ADE quiver varieties of Nakajima 
are shown to be subvarieties of the
variety of representations of a certain associative algebra
introduced by Lusztig. 
\end{abstract}

\maketitle

\section{Introduction}\label{Algebras}

\subsection{}\label{LusztigMap}
Motivated in part by G.~ Lusztig's work on quantum groups,
H.~ Nakajima introduced in \cite{Nakajima98}
remarkable algebraic varieties which we call
Nakajima's quiver varieties. 
In \cites{Lusztig98, Lusztig00, Lusztig03}
Lusztig produced another description of these varieties
as certain Grassmannians. The key element
of Lusztig's construction is an associative algebra attached to the quiver
(this is the Lusztig's algebra from the title)
and a map $\vartheta$ from Nakajima's varieties to
the variety of representations of this algebra.
In the $ADE$ case the map $\vartheta$ is proved to be finite and
a homeomorphism onto its image
\cite{Lusztig00}*{Theorem 4.7}.

The main purpose of the present paper is to prove 
(cf. Theorem \ref{ImmersionTheorem}) that if
the quiver is of $ADE$ type then 
the map $\vartheta$ is an algebraic isomorphism 
onto its image. This result is known in type 
$A$ (cf. \cite{Maffei98}) and was used in \cite{MirkovicVybornov03}
to construct a compactification of Nakajima's varieties
of type $A$. We plan to use techniques developed in 
the present paper to compactify Nakajima's varieties of arbitrary type.

\subsection{}\label{DefinitionOfAlgebras}
Although our motivation comes from the geometry of 
quiver varieties the main technical results of the
present paper are  purely algebraic and concern the preprojective
algebra of a graph. Let us recall the definitions.
Consider a finite graph $G$ with the vertex set $I$ and the
edge set $H$. 
We are mostly interested in graphs of finite and affine types, 
in other words in simply laced connected Dynkin 
or extended Dynkin diagrams, but at this point no restrictions on
$G$ are imposed. Now let us fix a field $k$ and a vector
$\lambda \in k^I$ and
consider the following three $k$-algebras associated to
$G$ (these will be our main objects of study):

\begin{itemize}
\item
the path algebra $kG$ which is the $k$-linear space 
spanned by the set of all paths in $G$ with multiplication 
(denoted by $\cdot$) given by concatenation of paths. We consider
$\lambda$ as an element of $kG$ under the natural identification of
the subalgebra of trivial paths with $k^I$. 
\item\label{DefinitionOfPi}
the deformed preprojective algebra $\Pi^{\lambda}$ which is the quotient
of $kG$ with respect to the two-sided ideal generated by the element
$(\theta - \lambda)$, where
\begin{equation}\label{DefinitionOfTheta}
\theta = \sum_{\substack{i \in I \\ a \in H}} 
\theta_i^a \medspace .
\end{equation}
Here $\theta_i^a \in kG$ is defined as follows.
If an edge $a$ has $i$ as an end-vertex then
$\theta_i^a$ is the path of length two from $i$ to $i$ along the
edge $a$ and back along the same edge, otherwise $\theta_i^a$ is zero.
\item\label{DefinitionOfL}
the Lusztig's algebra $L^{\lambda}$ which is isomorphic
to $kG$ as a $k$-linear space, but with the following
multiplication (cf. \cite{Lusztig03}*{2.1}):
\begin{equation}\label{LusztigProduct}
f \circ f' = f \cdot (\theta - \lambda ) \cdot f' \medspace , 
\end{equation}
where $\theta$ is as in \eqref{DefinitionOfTheta}. 
\end{itemize}

The structure of $\Pi^{\lambda}$ and 
$L^{\lambda}$ depends on the choice of
$\lambda$. For example
Crawley-Boevey and Holland have shown 
(see \cite{CrawleyBoeveyHolland98}*{Corollary 6.3}) that 
if $G$ is of finite type and
$\lambda$ does not lie on a Weyl chamber wall under the 
identification of $k^I$ with the $k$-span of the weight lattice 
then $\Pi^{\lambda}=0$ (zero algebra). Hence the most interesting 
case is $\lambda=0$ when we have the original
preprojective algebra $\Pi^0$. It is also explained in 
\cite{CrawleyBoeveyHolland98}
that in the affine case $\Pi^{\lambda}$
is related to the algebra of global functions
on a non-commutative deformation of a Kleinian singularity
(see below for a more precise statement).
The authors of \cite{CrawleyBoeveyHolland98} use slightly
different definition of the preprojective algebra, but it
is equivalent to ours for bipartite graphs
(cf. the discussion in \ref{QuiverAlgebras}).

\subsection{}
Our first theorem describes the Hilbert series
of the preprojective algebra of an arbitrary graph.
Let $\Pi^0 (n)$ denote the $n$-th graded component of
$\Pi^0$ with respect to the grading by path length. 
Given two vertices $i,j\in I$ let 
$\Pi^0_{ij} (n) = e_i \cdot \Pi^0 (n) \cdot e_j$
be the subspace of paths going from $j$ to $i$
and put $H_{ij} (n)=\dim \Pi^0_{ij} (n)$. The matrix Hilbert
series $H (t)$ of the preprojective algebra is defined as
\begin{equation}\nonumber
H_{ij}(t) = \sum_{n=0}^{\infty} H_{ij} (n) t^n \medspace .
\end{equation}

\begin{theorem}\label{HilbertTheorem}
Assume $\charac k =0$. Let $G$ be a connected graph and 
$A$ be the adjacency matrix of $G$.
\begin{enumerate}
\setcounter{enumi}{\value{equation}}
\item\label{IntroHilbertFinite}
If $A$ is not of ADET type then 
$H(t)=(1-At+t^2)^{-1}$. Moreover in this case the algebra $\Pi^0$
is Koszul.
\item\label{IntroHilbertAffine}
If $A$ is of ADET type then 
$H(t)=(1+Pt^h)(1-At+t^2)^{-1}$ where $h$ is the Coxeter number of $A$ and
$P$ is the permutation matrix corresponding to an involution of $I$.
\setcounter{equation}{\value{enumi}}
\end{enumerate}
\end{theorem}

Here $ADE$ refers to the simply laced Dynkin graphs and
$T_n$ is the following graph on $n$ vertices:
\begin{equation}\nonumber
\xygraph{
!{<0pt,0pt>;<15pt,0pt>:}
[] *\cir<2pt>{}
-@{-} [rr]
*\cir<2pt>{}
-@{{}.} [rr]
*\cir<2pt>{}
-@{-} [rr]
*\cir<2pt>{}
-@{-}@(ur,dr) []
*\cir<2pt>{}
}
\end{equation}
The Coxeter number of $T_n$ is $2n+1$.

Part (a) of the Theorem is known in many cases, see \cite{MartinezVilla96}. 
However our approach is completely different and we hope that it would
have other applications. Our observation is that one can think of the theory 
of preprojective algebras as of the invariant theory of the ``quantum plane''
with respect to the action of a ``finite subgroup of quantum $SL (2)$''.

Let us recall that McKay correspondence relates extended Dynkin
graphs and finite subgroups of $SL (2)$ and provides a 
very powerful approach to the study of affine quivers
(see e.g. 
\citelist{\cite{CrawleyBoeveyHolland98} \cite{Lusztig92} \cite{Lusztig98}}). 
Now one can replace $SL(2)$ with quantum $SL_q (2)$ and study
module categories over the category of representations of
$SL_q(2)$, the module categories playing the role of
subgroups in the quantum world.
It turns out that this quantum version of the McKay correspondence
allows one to obtain arbitrary graphs (cf. \cite{EtingofOstrik04}).
In the present paper we construct the preprojective algebra
using the quantum McKay correspondence and give a simple proof of
the above Theorem by reducing it to some standard facts about the quantum 
plane. In particular the special case of ADET graphs correspond to
the deformation parameter $q$ being a root of unity.

\subsection{}

Our second theorem describes
the zeroth Hochschild homology of the deformed
preprojective algebra in the finite and affine cases.

Let us recall that a prime number $p$ is called a \emph{bad prime}
for a simply laced Dynkin graph $G$ if
\begin{equation}\nonumber
\begin{split}
\text{$p=2$}        &\quad \text{ for $G$ of type $D_n$,} \\
\text{$p=2$ or $3$} &\quad \text{ for $G$ of type $E_6$ or $E_7$,} \\
\text{$p=2$, $3$, or $5$} &\quad \text{ for $G$ of type $E_8$.} 
\end{split}
\end{equation}
There are no bad primes for type $A$.
If $G$ is of affine type we say that $p$ is bad if
it is bad for the corresponding non-extended Dynkin graph.
A prime which is not bad is \emph{good}.
Good and bad primes play an important role in geometry
and representation theory of algebraic
groups. Their appearance in the following theorem is a
mystery for us.

\begin{theorem}\label{MainLemma}
Let $B$ be the image of the space of paths of length zero
under the canonical projection 
$kG \rightarrow \Pi^{\lambda}$.
\begin{enumerate}
\setcounter{enumi}{\value{equation}}
\item\label{LemmaFinite} 
Assume $G$ is of finite $ADE$ type and $\charac k$ is good for $G$.
Then
$\Pi^{\lambda} = [\Pi^{\lambda} , \Pi^{\lambda}] + B$.
\item\label{LemmaAffine}
Assume $G$ is of affine $\widehat{ADE}$-type 
and $\charac k$ is good for $G$.
Then
$\Pi^{\lambda} = [\Pi^{\lambda} , \Pi^{\lambda}] + 
\Pi^{\lambda}_{i_0 i_0} + B$,
where $\Pi^{\lambda}_{i_0 i_0} = 
e_{i_0} \cdot \Pi^{\lambda} \cdot e_{i_0}$ 
is the subalgebra
of $\Pi^{\lambda}$ consisting of all paths beginning and
ending at a chosen extending vertex $i_0 \in I$.
\end{enumerate}
\end{theorem}

Recall that if $R$ is an associative $k$-algebra then
$[R , R]$ denotes the $k$-span of the set
$\{ ab-ba \medspace \vert \medspace a , \medspace b \in R  \}$
and that a vertex of an affine Dynkin graph is called
extending if by removing it one obtains a connected
Dynkin graph of finite type.

The algebra $\Pi^{\lambda}_{i_0 i_0}$ 
which appears in the affine part
of the above theorem is very important. As explained in 
\cite{CrawleyBoeveyHolland98} it should 
be thought of as the algebra of global functions
on a non-commutative deformation of a Kleinian singularity.
If $\charac k = 0$ statement \ref{LemmaAffine} can also be
derived using McKay correspondence
(cf. \cite{CrawleyBoeveyEtingofGinzburg05}*{8.6}).

\subsection{}\label{QuiverAlgebras}
When $\lambda = 0$ the definition of the preprojective algebra 
contained in \ref{DefinitionOfAlgebras} is the original
one given by I. M. Gelfand and V. A. Ponomarev and
generalized by V. Dlab and C. M. Ringel. 
There is another definition in the literature
which uses quivers (oriented graphs) instead of non-oriented graphs.
Let us describe that construction. We fix an orientation $\Omega$ of
$G$ and denote by $Q$ the resultant quiver. We 
continue to denote by $H$ the set of edges (but now
each edge is, in fact, an arrow). Let $\overline{Q}$ be the double 
of $Q$, i.e. for each arrow $a \in H$ there is an extra arrow $a^*$
of $\overline{Q}$ connecting the same pair of vertices but going 
in the opposite direction. 

Let $k\overline{Q}$ be the path algebra of 
the oriented graph $\overline{Q}$ (if the graph $G$ has no self-loops
then it is canonically isomorphic to the path algebra $kG$)
and let $\theta_{\Omega}$ be the following element of
$k\overline{Q}$ (cf. \eqref{DefinitionOfTheta}):
\begin{equation}\label{DefinitionOfThetaOmega}
\theta_{\Omega} = \sum_{a \in H} (a \cdot a^* -a^* \cdot a) 
\medspace .
\end{equation}
Then the deformed preprojective algebra 
$\Pi^{\lambda}_{\Omega}$ of the \emph{quiver}
$Q$ is the quotient of $k\overline{Q}$ with respect to
the two-sided ideal generated by the element 
$(\theta_{\Omega} - \lambda)$ and
Lusztig's algebra $L^{\lambda}_{\Omega}$ is defined to be
the vector space $k \overline{Q}$ equipped with 
the multiplication 
\begin{equation}\nonumber
f \circ f' = f \cdot (\theta_{\Omega} - \lambda) \cdot f'
\medspace .
\end{equation}

It is easy to see that the two definitions of the preprojective
algebra are equivalent (i.e. $\Pi^{\lambda}$ and $\Pi^{\lambda}_{\Omega}$
are isomorphic) for bipartite graphs (in particular for all finite 
and affine graphs except $\Hat{A}_{2n}$). For arbitrary graphs
the algebras are different but our results remain valid. Namely,
the Hilbert series of $\Pi^{\lambda}_{\Omega}$ is given by
Theorem \ref{HilbertTheorem} (with $A$ being the adjacency
matrix of $\overline{Q}$) and its zeroth Hochschild homology 
is described by Theorem \ref{MainLemma}. We will indicate the required
modifications in the proofs. 

The crucial difference 
between non-oriented graphs and doubled quivers is not the
orientation and the related signs in the element $\theta$ but
rather the doubling of self-loops, which changes dramatically 
the structure of the preprojective algebra. For example,
if $G$ is of type $T_1$ (one vertex, one edge graph) then
$\Pi^{0} \simeq k[x]/(x^2)$ while 
$\Pi_{\Omega}^{0} \simeq k[x,y]$. 
Numerically each self-loop
of $G$ contributes 1 to the corresponding diagonal coefficient of
the adjacency matrix $A$ (which appears in Theorem \ref{HilbertTheorem})
but it adds 2 to the coefficient of the adjacency matrix of $\overline{Q}$.
Hence type $T_1$ is finite for $G$ (i.e. $(2I - A)$ is positive-definite)
but it is affine for $\overline{Q}$ (i.e. $(2I - A)$ is indefinite). The
Hilbert series of $\Pi^{0}_{\Omega}$ is given by the expression in
\eqref{IntroHilbertFinite} if $G$ is of $ADE$-type, and by
\eqref{IntroHilbertAffine} if it is not.

We decided to use the definition of the preprojective
algebra given in \ref{DefinitionOfAlgebras} because it works
more naturally in the context of the quantum McKay correspondence
(particularly if the graph has self-loops). However both
$\Pi^{\lambda}$ and $\Pi^{\lambda}_{\Omega}$ can be obtained
using module categories and we will indicate the required adjustments.

\subsection{}
The paper is organized as follows. 
In Section \ref{McKaySection} we describe a construction of
the preprojective algebra via quantum McKay correspondence
and prove Theorem \ref{HilbertTheorem}. 
An important example of star-shaped graphs is discussed in
Section \ref{StarSection}.
Section \ref{CalculationSection} contains the proof of 
Theorem \ref{MainLemma}.
Finally in Section \ref{VarietiesSection} we prove that Lusztig's
map $\vartheta$ (cf. \ref{LusztigMap}) is an isomorphism 
onto its image and describe applications to geometry
and Poisson structure of quiver varieties.

\subsection{} Recently another construction of 
quiver varieties appeared in 
\cite{FrenkelKhovanovSchiffmann03}. It uses 
so-called {\em skew-zigzag} algebra quadratic dual to the
preprojective algebra.
It would be very interesting to understand precisely the 
relationship (perhaps provided by Koszul duality) 
of the construction of 
\cite{FrenkelKhovanovSchiffmann03}
to Lusztig's construction and the present paper.

\subsection*{Acknowledgments}

We have greatly benefited from Igor Frenkel's inspiration 
and encouragement and  George Lusztig's course at MIT 
in the fall of 2002. We are grateful to P. Etingof for 
useful discussions and for careful reading of a previous version of 
the paper. Numerous comments and suggestions of the journal referee
are also very much appreciated.

A.~M.~ and M.~V. would like to thank the IH\' ES for 
their hospitality and support. The research of A.~M.~ and M.~V.~ was 
partially supported by the NSF Postdoctoral Research Fellowships. 
The research of  V.~O.~ was partially supported 
by the NSF grants DMS-0098830 and DMS-0111298.

\section{Quantum McKay correspondence and the
preprojective algebra}\label{McKaySection}

In this section we assume that $k$ is an algebraically
closed field of characteristic zero.

\subsection{Quantum McKay Correspondence}
As explained in the Introduction our goal is to obtain (and to exploit)
a construction of the preprojective algebra of an arbitrary graph via quantum
McKay correspondence. The original McKay observation was that given a 
finite subgroup $\Gamma \neq \{ \pm 1\}$ of $SL (2)$ 
one can associate to it a simply laced affine Dynkin graph as follows. Let
$V$ be the fundamental $2$-$\dim$ representation of $SL(2)$.
Then the McKay graph has simple representations of $\Gamma$ as vertices and
two representations $\rho$ and $\rho'$ are connected by 
$\dim \Hom_{\Gamma} (\rho', \rho \otimes V)$ edges. McKay correspondence
allows one to simplify many constructions and proofs in the particular
case of an affine graph. For example a construction of the preprojective
algebra of an affine quiver in terms of the invariant theory 
of $\Gamma$ is given in \cite{CrawleyBoeveyHolland98}*{Section 3}
and \cite{Lusztig98}*{Section 6}. 
Let us note that McKay theory relies only on the following facts: 
\begin{itemize}
\item 
finite dimensional representations of $SL(2)$ 
form a semisimple tensor category $\mathcal{C}$;
\item
representations of $\Gamma$ form a semisimple module category
$\mathcal{D}$ over $\mathcal{C}$ (i.e. there is a bifunctor
$\mathcal{C} \times \mathcal{D} \rightarrow \mathcal{D}$ satisfying
natural axioms, cf. \cite{Ostrik03}*{Section 2.3});
\item
there is a self-dual object $V \in \mathcal{C}$ (the fundamental
representation) such that
the composition
\begin{equation}\label{Trace}
\mathbf{1} \xrightarrow{\coev_V} V \otimes V^* 
\xrightarrow{\phi \otimes \phi^{-1}} V^* \otimes V
\xrightarrow{\ev_V} \mathbf{1} 
\end{equation}
is equal to $-2 \id_{\mathbf{1}}$ for any choice of an isomorphism
$\phi: V \rightarrow V^*$.  
\end{itemize}

Following \cite{EtingofOstrik04} we generalize this setup as
follows. Let $q \in k \setminus \{ 0 \}$
and consider the following tensor category $\mathcal{C}_q$: 
for $q=\pm 1$ or $q$ not equal to a root of unity $\mathcal{C}_q$ 
is the category of representations of the quantum 
group $SL_q(2)$ while for $q$ equal to a root of unity  
$\mathcal{C}_q$ is the semisimple
subquotient of the category of representations of $SL_q(2)$.

Let us recall the structure of $\mathcal{C}_q$ 
(cf. \cite{Kassel95}). First assume that $q$ is not a root of unity.
Then for each $s \in \mathbb{N}$ there is a simple
object $L_s$ of $\mathcal{C}_q$ (the deformation of the $(s+1)$-dimensional
representation of $SL(2)$), the set 
$\{ L_s \}_{s=0}^{\infty}$ is the complete set of simple objects, and
the tensor product decomposition is given by
\begin{equation}\label{TensorGeneric}
L_r \otimes L_s \simeq 
\bigoplus_{\substack{t = |r-s| \\ t \equiv r+s \mod 2}}^{r+s} L_t
\medspace .
\end{equation}

If $q$ is a root of unity we put
\begin{equation}\label{CoxeterQ}
h(q) = 
\begin{cases}
N      &\text{ if $q$ is a root of unity of an odd order $N$,} \\
\frac{N}{2} &\text{ if $q$ is a root of unity of an even order $N$.} 
\end{cases} 
\end{equation}
Then for each integer $s=0, \ldots , h(q)-2$ there is a simple
object $L_s$ of $\mathcal{C}_q$ (the deformation of the $(s+1)$-dimensional
representation of $SL(2)$), the set 
$\{ L_s \}_{s=0}^{h(q)-2}$ is the complete set of simple objects, and
the tensor product decomposition is given by
\begin{equation}\label{TensorRoot}
L_r \otimes L_s \simeq 
\begin{cases}
\displaystyle
\bigoplus_{\substack{t = |r-s| \\ t \equiv r+s \mod 2}}^{r+s} 
L_t \quad &\text{ if } r+s < h(q)-1 \medspace ; \\
\\
\displaystyle
\bigoplus_{\substack{t = |r-s| \\ t \equiv r+s \mod 2}}^{2h(q)-4-r-s} 
L_t \quad &\text{ if } r+s \geq h(q)-1 \medspace .
\end{cases}
\end{equation}

Let $V=L_1$ be the fundamental representation of $SL_q(2)$.
The main result of \cite{EtingofOstrik04} is a classification of
semisimple module categories over $\mathcal{C}_q$ with finitely
many isomorphism classes of simple objects. 
Let $\mathcal{D}$ be such a category and $I$ be the set of 
isomorphism classes of simple objects of $\mathcal{D}$. 
Then as an abelian category
$\mathcal{D}$ is equivalent to the category of $I$-graded vector spaces
which we denote by $\mathcal{M}_I$. 
Let $\Fun(\mathcal{M}_I,\mathcal{M}_I) \cong 
\mathcal{M}_{I\times I}$ 
be the category of additive functors from
$\mathcal{M}_I$ to itself. The category 
$\Fun(\mathcal{M}_I,\mathcal{M}_I)$ has an obvious structure
of monoidal category induced by the composition of functors and
the structure of a $\mathcal{C}$-module category 
on $\mathcal{M}_I$ is just a tensor functor 
$F:\mathcal{C}_q \to \Fun(\mathcal{M}_I,\mathcal{M}_I)$.
According to \cite{EtingofOstrik04}, if $q$ is not a root of unity then 
such functors are classified by the following data:
\begin{itemize}
\item
a collection of finite dimensional 
vector spaces $\{ V_{ij}\}_{i,j\in I}$;
\item
a collection of nondegenerate bilinear forms 
$E_{ij}: V_{ij}\otimes V_{ji}\to k$
satisfying the following conditions 
\begin{equation}\label{EEEquation}
\sum_{j\in I} \Tr (E_{ij}(E^t_{ji})^{-1})=-q-q^{-1} 
\end{equation}
for each $i \in I$.
\end{itemize}
Here $\bigoplus_{i,j\in I}V_{ij} = F(V)$ is the object of 
$\mathcal{M}_{I\times I} \cong 
\Fun (\mathcal{M}_I,\mathcal{M}_I)$ 
corresponding to $V \in \mathcal{C}_q$ and the forms $E_{ij}$ are
induced by an isomorphism $\phi: V \rightarrow V^*$. 
If $q$ is a root of unity there is an extra condition
due to the fact that $\mathcal{C}_q$ is a quotient of the full 
tensor category generated by $V$. To avoid technical details we don't
spell out the extra condition here and refer the reader to
\cite{EtingofOstrik04} instead. 

Note that since $E_{ij}$ is nondegenerate we have 
$\dim V_{ij} = \dim V_{ji}$.
Let $A$ be the symmetric $I \times I$ matrix given by
$A_{ij} = \dim V_{ij}$. As in the original McKay correspondence
we think of $A$ as the adjacency matrix of a graph with vertex set $I$. 
Let us note that graphs equipped with a collection of bilinear forms
assigned to edges are called modulated graph and they 
have been studied before in connection with non simply
laced root systems (cf. \cite{DlabRingel76}). So the quantum
McKay correspondence produces modulated graphs satisfying condition
\eqref{EEEquation}.

The crucial fact is that by choosing an appropriate value of $q$ and
an appropriate  module category $\mathcal{D}$ we can
obtain an arbitrary connected graph in this way. More precisely,
according to \cite{EtingofOstrik04}, given a connected
graph (i.e. dimensions of $V_{ij}$'s)
the equation \eqref{EEEquation} has a solution (not necessarily unique)
for some choice of $q$
(in other words there exists a modulation of the graph satisfying
\eqref{EEEquation}). 
We shall describe a solution to \eqref{EEEquation} in the course of 
the proof of Lemma \ref{PiELemma}.

It is shown in \cite{EtingofOstrik04} that
$q$ is a root of unity if and only if the McKay graph is a 
Dynkin graph of ADET type with Coxeter number $h(q)$ 
(cf. \eqref{CoxeterQ}). 

\subsection{Symmetric algebra}
Let $T$ be the free algebra in the category $\mathcal{C}_q$ 
generated by $V$. Under the functor 
$F: \mathcal{C}_q \rightarrow \mathcal{M}_{I \times I}$ it maps to
the path algebra of the McKay graph, i.e. to the algebra generated by 
the set of orthogonal idempotents $\{ e_i \}_{i \in I}$ and
the space of generators of degree one $\oplus_{i, j \in I} V_{ij}$.

Now consider the quotient $S$ of $T$ with respect to 
the two sided ideal generated
by the image of $\mathbf{1}$ under the map
$\mathbf{1} \xrightarrow{\coev_V} V \otimes V^* 
\xrightarrow{\id_V \otimes \phi^{-1} } V \otimes V$.
If $q=1$ (i.e. in the classical situation)
$S$ is just the algebra of polynomials in two
commuting variables. In general
$S$ is called the $q$-symmetric algebra
or the algebra of functions on the quantum plane.
Applying the functor $F$ one gets an algebra $\tilde{\Pi}^E =F(S)$ 
which is the quotient
of the path algebra with respect to the two sided ideal generated by 
$F (\mathbf{1} \subset V \otimes V)$. By the definition of $F$ we have
\begin{equation}\nonumber
F (\mathbf{1} \subset V \otimes V) = \sum_{i \in I} \theta^E_i 
\medspace ,
\end{equation}
where $\theta^E_i$ is defined as follows.
Fix a pair of bases $x_1, \ldots , x_{A_{ij}}$ in $V_{ij}$ and 
$y_1 , \ldots , y_{A_{ij}}$ in $V_{ji}$ dual with respect to the form
$E_{ji}$ and let $\theta^E_{ij}= \sum_{s = 1}^{A_{ij}} y_s \cdot x_s$. Then
the element $\theta^E_{ij}$ does not depend on the choice of the dual bases
and we put $\theta^E_i = \sum_{j\in I} \theta^E_{ij}$.
The algebra $\tilde{\Pi}^E$ is called the preprojective algebra of the 
modulated graph (cf. \cite{DlabRingel80}). 

\begin{lemma}\label{PiELemma}
Let $G$ be an arbitrary connected graph and $\Pi^0$ be
its preprojective algebra defined in \ref{DefinitionOfPi}. 
Then there exists a value of $q$ and a $\mathcal{C}_q$-module category
(in other words, a modulation $\{ E_{ij} \}$ of $G$) 
such that $\Pi^0 = \tilde{\Pi}^E$.
\end{lemma}
The lemma means that the preprojective algebra of an arbitrary
connected graph can be realized via quantum McKay correspondence, or
in other words, that the algebra of functions on the quantum
plane is the ``universal preprojective algebra''.
\begin{proof}
Let $A$ be the adjacency matrix of $G$.  
It is irreducible since $G$ is connected.
Let $\lambda$ be the Frobenius-Perron eigenvalue of $A$ with eigenvector
$\{ r_i \}_{i \in I}$ (in particular $r_i \neq 0$ for any $i \in I$)
and choose $q$ so that $\lambda = -q-q^{-1}$.
Now we can define the module category. Given two vertices 
$i$, $j \in I$ let $V_{ij}$ be the linear space generated 
by the set of edges between $i$ and $j$ and define
the bilinear forms
$E_{ij}: V_{ij} \otimes V_{ji} \rightarrow k$ as follows:
\begin{equation}\label{EIJFormulas}
E_{ij} (a , b) = 
\begin{cases}
r_j  &\text{ if $a=b$ ,} \\
0    &\text{ if $a \neq b$ ,}
\end{cases}
\end{equation}
where $a$ and $b$ are two edges between $i$ and $j$ considered as basic
elements of $V_{ij}$ and $V_{ji}$ respectively.

We claim that the forms $E_{ij}$ satisfy the condition \eqref{EEEquation}.
Indeed,
\begin{multline}\nonumber
\sum_{j\in I} \Tr (E_{ij}(E^t_{ji})^{-1}) = 
\sum_{j\in I} \dim V_{ij} \Bigl( \frac{r_j}{r_i} \Bigr)
= \\ =
\sum_{j\in I} A_{ij} \frac{r_j}{r_i} 
= \lambda \frac{r_i}{r_i}
= - q - q^{-1} 
\end{multline}
since $\{ r_i \}_{i \in I}$ is an eigenvector of $A$ with eigenvalue 
$\lambda$. 

Now with this choice of bilinear forms we have 
\begin{equation}\nonumber
\theta^E_i = \frac{1}{r_i} \sum_{a \in H} \theta_i^a \medspace ,
\end{equation}
where $\theta^a_i$ is as in \eqref{DefinitionOfTheta}.
It follows that $\tilde{\Pi}^E = \Pi^0$.
\end{proof}

\subsection{Hilbert Series.}

The prescription $\deg e_i = 0$ and $\deg x = 1$ for 
$x\in V_{ij}$ defines a grading on the algebra 
$\tilde{\Pi}^E$. Let $\tilde{\Pi}^E (n)$
denote the $n$-th graded component. Given $i,j\in I$ let 
$\tilde{\Pi}^E_{ij} (n) = e_i \tilde{\Pi}^E (n) e_j$
and let $H_{ij} (n)=\dim \tilde{\Pi}^E_{ij} (n)$. 
In this subsection we calculate the matrix Hilbert
series $H(t)=\sum_{n=0}^\infty H_{ij}(n) t^n$ in terms of the
adjacency matrix $A$.

\begin{theorem}\label{HilbertTensorThm}
Let $H(t)$ be the matrix Hilbert series of $\Tilde{\Pi}^E$
defined above.
\begin{enumerate}
\setcounter{enumi}{\value{equation}}
\item\label{HilbertKM}
Assume that the adjacency matrix $A$ is not of ADET type. Then 
$H(t)=(1-At+t^2)^{-1}$. Moreover in this case the algebra 
$\Tilde{\Pi}^E$ is Koszul.

\item\label{HilbertFinite}
Assume that the adjacency matrix $A$ is of $ADET$ type. Then 
$H(t)=(1+Pt^h)(1-At+t^2)^{-1}$, where $h$ is the Coxeter number of $A$ and
$P$ is the permutation matrix corresponding to some involution of $I$.
\setcounter{equation}{\value{enumi}}
\end{enumerate}
\end{theorem}

Thanks to Lemma \ref{PiELemma} the above Theorem is a generalization of 
Theorem \ref{HilbertTheorem} and in particular implies the latter. 
A minor technical point is that in Theorem \ref{HilbertTheorem}
we don't assume $k$ to be algebraically closed. However it is clear
that $H(t)$ does not change under base field extensions.

\noindent\textbf{Remark.} It is easy to show that 
the Hilbert series of the quadratic dual 
of $\tilde{\Pi}^E$ is equal to $1 + At + t^2$
(except for several trivial graphs). 
It is also known that the Hilbert polynomial
of the Yoneda algebra of $\Tilde{\Pi}^E$ is given by 
$1 + At + t^2$ (for example it follows from the last
part of the proof below). Therefore the algebra $\tilde{\Pi}^E$
is Koszul if and only if its Hilbert series
is given by $(1-At+t^2)^{-1}$.

\begin{proof} 
The structure of the $q$-symmetric algebra $S$  
is well known, see \cite{Kassel95}. In the case 
\eqref{HilbertKM} (so $q$ is not a root of unity) 
we have $S (n) = L_n$ for any $n\in \mathbb{N}$, where 
$S (n)$ is the $n$-th graded component of $S$. 
Now it follows from \eqref{TensorGeneric} that
\begin{equation}\nonumber
L_1\otimes L_n \simeq L_{n-1}\oplus L_{n+1}
\qquad \text{for $n>0$,}
\end{equation} 
and applying the functor $F$ we get a recursion
\begin{equation}\nonumber
tAH(t) = H(t) + t^2H(t) \medspace ,
\end{equation}
which implies the formula in \eqref{HilbertKM}.

In the case \eqref{HilbertFinite} (so $q$ is a root of 1)
we have $S (n)=L_n$ for $0 \le n \le h-2$ and
$S (n)=0$ for $n\ge h-1$. Let 
\begin{equation}\nonumber
\Hat{S} = S \ominus t^h (L_{h-2} \otimes S) \oplus  
t^{2h} (L_{h-2} \otimes L_{h-2} \otimes S) \ominus \ldots \medspace ,
\end{equation}
where powers of $t$ denote grading shifts. Then
\eqref{TensorRoot} implies the following recursion 
\begin{equation}\nonumber
L_1 \otimes \Hat{S}(n) \simeq \Hat{S}(n-1) \oplus \Hat{S}(n+1)
\qquad \text{for $n>0$},
\end{equation} 
and applying the functor $F$ we obtain the formula in 
\eqref{HilbertFinite} with $P = F ( L_{h-2} )$. Note that  
$L_{h-2}\otimes L_{h-2} \simeq L_0 = \mathbf{1}$, hence $P^2 = \id$.

Now assume that $q$ is not a root of unity or $q=\pm 1$. Then the trivial 
left module 
$\mathbf{1}=L_0$ over the algebra $S$ has the following resolution:
\begin{equation}\nonumber
\mathbf{1} \rightarrow S \rightarrow 
S \otimes L_1 \rightarrow S \rightarrow \mathbf{0}
\end{equation}
(this is just the $q-$version of the usual resolution for the algebra of 
polynomials in two commuting variables).
The image of this resolution under the functor $F$ provides a
Koszul resolution for the algebra $\tilde{\Pi}^E$. 
\end{proof}

\noindent\textbf{Remark.} The involution $P$ in \eqref{HilbertFinite} is an 
automorphism of the underlying graph and is explicitly known in all cases:
for types $A_n, D_{2n+1}, E_6$ it is the unique nontrivial involution 
while for $D_{2n}$, $T_n$, $E_7$, $E_8$ it is the identity map.

\subsection{} In the $ADET$ case the preprojective algebra is finite
dimensional and we have the following result.

\begin{corollary}\label{ADETCorollary}
If the graph $G$ is of $ADET$-type then
\begin{equation}\nonumber
\dim \Pi^{0}=\frac{h(h+1)r}{6} \medspace .
\end{equation}
where $h$ and $r$ are the Coxeter number 
and the rank respectively. 
\end{corollary}

\begin{proof}
One has to put $t=1$ in \eqref{HilbertFinite} and
then calculate the sum of the matrix elements of the resultant matrix. 
In other words 
$\dim \Pi^0 = \Tr(H(1)U)$, where the matrix $U$ is given by
$U_{ij}=1$ for all $i,j\in I$. Observe that $PU=U$ and hence
$\dim \Pi^0 = 2\Tr((2-A)^{-1}U)$. Now assume we are in $ADE$ situation.
Then $(2-A)$ is the Cartan
matrix of the corresponding root system and hence $(2-A)^{-1}_{ij}=
(\omega_i,\omega_j)$ where $\omega_i$ are the fundamental weights and the
scalar product is normalized by the condition $(\alpha,\alpha)=2$ for a
root $\alpha$. 
Therefore $\dim \Pi^0=2\sum_{i,j\in I}(\omega_i,\omega_j)=2(\rho,\rho)$, 
where
$\rho =\sum_{i\in I}\omega_i$. Now the result follows from the ``strange
formula'' of Freudenthal and de Vries.

If $G$ is of type $T_n$ one can think of it as of a ``folded'' graph
$\widetilde{G}$ of type $A_{2n}$. Then we have the natural surjective 
map of path algebras $\sigma: k\widetilde{G} \rightarrow kG$.
The map $\sigma$ respects the grading and
the relations \ref{DefinitionOfTheta}, and it is injective on the space 
of paths beginning at a fixed vertex of $\widetilde{G}$. 
Since each vertex of $G$ has two preimages in
$\widetilde{G}$ we have
$\dim \Pi^0 (\widetilde{G}) = 2 \dim \Pi^0 (G)$ which implies 
the Corollary (recall that the Coxeter numbers of $G$ and 
$\widetilde{G}$ are
the same while the rank of $\widetilde{G}$ is twice the rank of $G$).
\end{proof}

\subsection{Quiver algebras}As discussed in Section \ref{QuiverAlgebras},
Theorem \ref{HilbertTheorem} remains valid for the 
preprojective algebra $\Pi^{0}_{\Omega}$ defined using a double
quiver $\overline{Q}$ instead of a graph $G$ (with $A$ being the
adjacency matrix of $\overline{Q}$). In fact one can construct
the algebra $\Pi^{0}_{\Omega}$ using quantum McKay correspondence. Namely,
the following variant of Lemma \ref{PiELemma} holds true and implies the 
analogue of Theorem \ref{HilbertTheorem}.

\begin{lemma}
Let $Q$ be an arbitrary connected quiver. 
Then there exists a value of $q$ and a $\mathcal{C}_q$-module category 
such that $\Pi^0_{\Omega} = \tilde{\Pi}^E$.
\end{lemma}

\begin{proof}The proof is similar to the proof of Lemma \ref{PiELemma}.
To define the module category we put $V_{ij}$ to be the linear space
generated by arrows of $\overline{Q}$ going from $i$ to $j$ and define
the bilinear forms $E_{ij}$ as follows
(recall that $H$ is the set of arrows of $Q$):
\begin{equation}\nonumber
E_{ij} (a,b) = 
\begin{cases}
r_j  &\text{ if $a \in H$ and $b=a^*$ ,} \\
-r_j  &\text{ if $b \in H$ and $a=b^*$ ,} \\
0    &\text{ otherwise .}
\end{cases} 
\end{equation}
These data define a $\mathcal{C}_q$-module category 
provided $\{ r_i \}_{i \in I}$ is the
Frobenius-Perron eigenvector of the adjacency matrix $A$ with the eigenvalue
$\lambda = q+q^{-1}$. Moreover one has $\Pi^0_{\Omega} = \tilde{\Pi}^E$. 
The calculations are
similar to the ones in the proof of Lemma \ref{PiELemma}.
\end{proof}

If $A$ is affine then $\lambda = 2$ and $q=1$. 
Thus we are in the classical McKay situation. 
In this case the above construction reduces to
McKay description of the preprojective algebra
given in  \cite{CrawleyBoeveyHolland98}*{Section 3}
and \cite{Lusztig98}*{Section 6}. 
Our choice of the bilinear forms $E_{ij}$ is exactly
the same as in \emph{loc. cit.} 
(where $r_i$ is the dimension of the
irreducible representation $\rho_i$ of $\Gamma$).

\section{Star graphs}\label{StarSection} 

In this section we study star-shaped graphs. 
We assume $\charac k = 0$.

\subsection{Spherical subalgebra}
A star graph looks as follows:
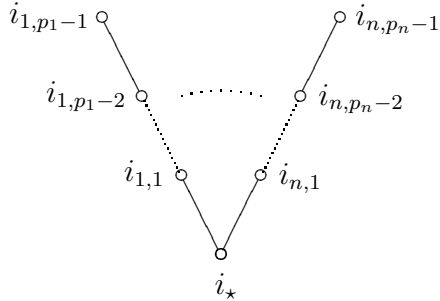
\begin{equation}\nonumber
\xygraph{
!{<0pt,0pt>;<15pt,0pt>:}
[] *\cir<2pt>{}
-@{-{}} [uul]
*\cir<2pt>{}
!{\save -<14pt,0pt>*\txt{$i_{1,1}$} 
\restore}
-@{{}.{}} [uul]
*\cir<2pt>{}
!{\save -<21pt,0pt>*\txt{$i_{1,p_1 - 2}$} 
\restore}
-@{-{}} [uul]
*\cir<2pt>{}
!{\save -<20pt,0pt>*\txt{$i_{1,p_1 - 1}$} 
\restore}
[ddddddrrr]
*\cir<2pt>{}
!{\save +<2pt,-12pt>*\txt{$i_{\star}$} 
\restore}
-@{-{}} [uur]
*\cir<2pt>{}
!{\save +<15pt,-1pt>*\txt{$i_{n,1}$} 
\restore}
-@{{}.{}} [uur]
*\cir<2pt>{}
!{\save +<23pt,-1pt>*\txt{$i_{n,p_n - 2}$} 
\restore}
-@{-{}} [uur]
*\cir<2pt>{}
!{\save +<22pt,-1pt>*\txt{$i_{n,p_n - 1}$} 
\restore}
[ddddddlll]
[uuuul]
-@/^2pt/@{{}.{}} [rr]
}
\end{equation}
It has $n$ rays of lengths
$p_1, \dots, p_n$ attached to the central vertex $i_\star$. 
Graphs of this shape are important for two reasons.
First they include all finite ADE types as well as some affine ones
($\Hat{D}_4$, $\Hat{E}_6$, $\Hat{E}_7$, and $\Hat{E}_8$). Second 
they play an important role in classical linear algebra,
such as various forms of Deligne-Simpson problem, saturation conjecture for
Littlewood-Richardson coefficients, etc.. We would like to state this 
connection to matrix problems in the form of the following Proposition
concerning the spherical subalgebra 
$\Pi_{i_{\star} i_{\star}}^0 = e_{i_\star} \cdot \Pi^0 \cdot e_{i_\star}$
of the preprojective algebra $\Pi^0$ consisting of paths 
beginning and ending at the central vertex.

\begin{proposition}\label{StarLemma} 
The unital algebra $\Pi_{i_{\star} i_{\star}}^0$ 
has the following presentation over $k$:
\begin{equation}\nonumber
\Pi_{i_{\star} i_{\star}}^0 =
< \ x_r,\ \ r=1, \dots ,n \ | \ x_1^{p_1}=\dots = x_n^{p_n}=
x_1+ \dots + x_n= 0 \ > \medspace .
\end{equation}
The presentation is given by sending $x_r$ to the loop of length one 
along the $r$-th ray of the graph
(i.e. $x_r = [i_\star i_{r,1} i_\star ]$).
\end{proposition}
\begin{proof}
Let $R_s$ be the graph
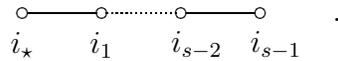
\begin{equation}\nonumber
\xygraph{
!{<0pt,0pt>;<15pt,0pt>:}
[] 
*\cir<2pt>{}
!{\save -<0pt,12pt>*\txt{$i_\star$} \restore}
-@{-{}} [rr]
*\cir<2pt>{}
!{\save -<0pt,12pt>*\txt{$i_1$} \restore}
-@{.{}} [rr]
*\cir<2pt>{}
!{\save -<-6pt,12pt>*\txt{$i_{s-2}$} \restore}
-@{-{}} [rr]
*\cir<2pt>{}
!{\save -<-6pt,12pt>*\txt{$i_{s - 1}$} \restore}
} \quad .
\end{equation}
One should think of $R_s$ as of a ray of the star graph.
Consider the algebra $\Phi_s$ defined as
\begin{equation}\nonumber
\Phi_s = k R_s \biggl/ 
\Bigl( {\textstyle \sum_{m=1}^{s-1}} \theta_{i_m} \Bigr)
\medspace , 
\end{equation}
where $\theta_{i_m}$ is as in \eqref{DefinitionOfTheta}.
The algebra $\Phi_s$ is almost the preprojective algebra
of the graph $R_s$ except we don't impose the relation 
$\theta_{i_{\star}}$. Let 
$\Psi_s = e_{i_\star} \cdot \Phi_s \cdot e_{i_\star}$
be the subalgebra of paths beginning and ending at the star vertex. 
An easy induction on $s$ shows that
\begin{equation}\nonumber
\Psi_s = k[x] \Bigl/ (x^s) \medspace 
\end{equation}
with $x = [i_\star i_1 i_\star ]$.
Now the Proposition follows from the fact that
for the original star graph
\begin{multline}\nonumber
\Pi^0_{i_{\star} i_{\star}} =
\Psi_{p_1} \ast \ldots \ast \Psi_{p_n} \Bigl/ ( \theta_{i_\star} ) 
= \\ =
\bigl( k[x_1]/(x_1^{p_1}) \bigr) \ast \ldots \ast 
\bigl( k[x_n]/(x_n^{p_n}) \bigr) \Bigl/ 
( x_1 + \ldots + x_n)  
\medspace ,
\end{multline}
where $\ast$ denotes the free product of algebras.
\end{proof}

\subsection{Hilbert Series.} 
We saw that the calculation of the Hilbert series of a
preprojective algebra reduces to the inversion 
of the matrix $1-At+t^2$. 
In the case of a star-shaped graph Lusztig and Tits 
\cite{LusztigTits92} derived
the following expression for the ``central'' 
coefficient of the inverse:
\begin{equation}\nonumber
\bigl((1-At+t^2)^{-1}\bigr)_{i_\star i_\star} =
\Bigl( 1 + t^{2}-t \sum_{s=1}^n 
\frac{t^{p_s -1} -t^{-p_s+1}}{t^{p_s}-t^{-p_s}} \Bigr)^{-1} \ .
\end{equation}
Using Theorem \ref{HilbertTheorem} and the above formula 
we obtain the following result (note that $\deg x_r = 2$ in the 
standard path algebra grading).

\begin{proposition}Let $h(t)$ be the Hilbert series of the subalgebra
$\Pi^0_{i_{\star} i_{\star}}$ with respect to the grading defined by 
$\deg(x_r)=1$.
\begin{enumerate}
\setcounter{enumi}{\value{equation}}
\item
If the star graph is not of ADE type then 
\begin{equation}\nonumber
h(t)= \Bigl( 1+t-\sum_{s=1}^n\frac{t-t^{p_s}}{1-t^{p_s}}
\Bigr)^{-1} \ .
\end{equation}
\item\label{HilbertADE}
If the star graph is of DE type then
\begin{equation}\nonumber
h(t)=
\Bigl( 1+t-\sum_{s=1}^n\frac{t-t^{p_s}}{1-t^{p_s}}
\Bigr)^{-1}(1+t^{h/2}) \ .
\end{equation}
\setcounter{equation}{\value{enumi}}
\end{enumerate} 
\end{proposition}
As an example, for $E_8$ we obtain
\begin{multline}\nonumber
h(t)= \Bigl( 1+t-\frac{t-t^2}{1-t^2}-\frac{t-t^3}{1-t^3}-
\frac{t-t^5}{1-t^5} \Bigr)^{-1}
(1+t^{15})= \\
1+2t+3t^2+4t^3+5t^4+6t^5+6t^6+6t^7+6t^8+6t^9+5t^{10}+4t^{11}+3t^{12}+2t^{13}
+t^{14}
\end{multline}

\subsection{Kleinian groups}
Finally we would like to explain the connection 
between star shaped graphs
and the McKay correspondence. 
Let us recall that finite subgroups of
$SL (2)$ are in one-to-one
correspondence with classical Dynkin graphs of ADE type.
Let $a$, $b$, and $c$, be the lengths of the rays of the graph
corresponding to a subgroup $\Gamma \subset SL(2)$. 
In the $A_{2n-1}$-case we put $a=b=n$ and $c=1$ and we exclude
$A_{2n}$-case to simplify discussion.
Then $-I \in \Gamma$ and we have the
following presentation (cf. \cite{Coxeter91}):
\begin{equation}\nonumber
\Gamma / \{ \pm I \} = <\ X , Y, Z \ | \ X^a=Y^b=Z^c=XYZ=1 \ > 
\medspace ,
\end{equation}
which is a multiplicative analogue of the presentation of
the spherical subagebra.
\begin{proposition}\label{StarLemma3}
The algebra $\Pi_{i_{\star} i_{\star}}^0$ has the following presentation
\begin{equation}\nonumber
\Pi_{i_{\star} i_{\star}}^0 
= < \ x, y, z \ | \ x^a=y^b=z^c=x+y+z=0 \ > \medspace .
\end{equation}
Moreover $\dim \Pi^0_{i_{\star} i_{\star}} = | \Gamma |/2$.
\end{proposition}

\begin{proof}
The first statement is a special case of Proposition \ref{StarLemma}. 
The dimension of $\Pi^0_{i_{\star} i_{\star}}$ in the $DE$ case
is given by \ref{HilbertADE} with $t=1$ and in the $A_{2n-1}$ case 
can be easily found from the presentation. The answer is the same
in both cases: 
\begin{equation}\nonumber
\dim \Pi^0_{i_{\star} i_{\star}} = 
\frac{2}{\frac{1}{a} + \frac{1}{b} + \frac{1}{c} - 1} 
= | \Gamma |/2 \ . 
\end{equation}
\end{proof}
\noindent\textbf{Remark.} 
As it was explained to us by P. Etingof one can both prove and
generalize the equality 
$\dim \Pi^0_{i_{\star} i_{\star}} = | \Gamma |/2$
using the holonomy of the connection 
$(\frac{x}{\xi} + \frac{y}{\xi + 1} + \frac{z}{\xi - 1}) d \xi$
on $\mathbb{C} \setminus \{0, \pm 1\}$.

\section{Proof of Theorem \ref{MainLemma}}\label{CalculationSection}

\subsection{}

This section is devoted to the proof of Theorem \ref{MainLemma}
(the same argument works for preprojective algebras of double quivers --
cf. \ref{QuiverAlgebras}).
The proof consists of several reductions followed with case-by-case 
calculations.
First note that since $\Pi^{\lambda}$ is filtered by the path length
it is enough to prove the Theorem for the associated graded
algebra which is isomorphic to a quotient of $\Pi^{0}$. 
So we assume $\lambda =0$.

\subsection{}

Next we do a reduction on the support of (a class of) a path.
Let $i$ be a vertex of $G$, $\Pi^0_{ii}=e_i \cdot \Pi^0 \cdot e_i$ 
be the subalgebra
of $\Pi^{0}$ consisting of classes of paths beginning and
ending at $i$, $\Pi^{0}_{G \setminus i}$ be the
preprojective algebra of the graph obtained by removing 
the vertex $i$ and the adjacent edges from $G$, 
and $J_i$ be the two-sided ideal of $\Pi^{0}$ generated by $e_i$
(so $J_i$ consists of paths ``passing through $i$''). 
We claim that $J_i \subset \Pi^{0}_{ii} + [\Pi^{0} , \Pi^{0}]$.
Indeed any non-cyclic path belongs to $[\Pi^{0} , \Pi^{0}]$ and any
cyclic path passing through $i$ can be made to begin and end at $i$ 
after a cyclic permutation (i.e. after adding an element of 
$[\Pi^{0} , \Pi^{0}]$). 
Now we have a natural surjective map of algebras 
$\Pi^{0}_{G \setminus i} \rightarrow \Pi^{0}/J_i$ and hence
\eqref{LemmaAffine} follows from \eqref{LemmaFinite}
if we put $i$ to be an extending vertex
while \eqref{LemmaFinite} can be proven by induction on
the number of vertices once it is shown that
$\Pi^{0}_{ii} \subset [\Pi^{0} , \Pi^{0}] +B$ for
some vertex $i$.
The verification of this last inclusion is the following
case-by-case calculation.

\emph{$G$ is of type $A_n$.} We choose $i$ to be one of the two
end-points of the Dynkin graph. Then according to (a degenerate
version of) Proposition \ref{StarLemma} we have
\begin{equation}\nonumber
\Pi^{0}_{ii} = < \ x \medspace | \medspace x^n = x = 0 \ > = k e_i \subset B
\end{equation}  

\emph{$G$ is of type $D_n$ or $E_n$.} We choose $i$ to be 
an end-point of the Dynkin graph furthest from the three-valent 
point. Let $f \neq e_i$ be a path in $G$ beginning and ending
at $i$. We are going to prove that $f =0$ modulo the commutant and
the ideal \ref{DefinitionOfTheta}. In what follows 
all the calculations are done
modulo these subspaces. First of all we may assume that 
$f$ passes through the three-valent vertex (otherwise we are
in type $A$ situation). Therefore (modulo the commutant) we assume
that $f$ begins and ends at the three-valent vertex and passes through
the end-point. We will use the notation of Proposition
\ref{StarLemma3} with $(a,b,c) = (n-2,2,2)$ (resp. $(n-3,3,2)$)
in type $D_n$ (resp. $E_n$). 
So $f$ is a word in the alphabet $\{ x,y,z \}$
containing $x^{a-1}$. We have to prove that $f=0$ if 
\begin{equation}\nonumber
x^a = y^b = z^c = x + y + z = 0
\end{equation} 
and we are allowed to do cyclic 
permutations. In other words, $f$ is a \emph{cyclic} word
in $y$ and $z$ subject to the following relations:
\begin{gather}\label{GeneralYRelation}
y^b = 0 \\ \label{GeneralZRelation}
z^c = 0 \\ \label{GeneralYZRelation}
(y+z)^a = 0
\end{gather}
and we want to show that $f=0$ if $f = (y+z)^{a-1}g$ for
some $g \neq 1$ (otherwise we are again in type $A$ situation).
Moreover because of \eqref{GeneralYZRelation}, we can assume
that the word $g$ begins and ends with $z$. 

\subsection{$D_n$-case}
We have 
\begin{gather}\label{DnYRelation}
y^2 = 0 \\ \label{DnZRelation}
z^2 = 0 \\ \label{DnYZRelation}
(y+z)^{n-2} = 0
\end{gather}
which implies that
\begin{equation}\label{Binom}
(y+z)^k = \underbrace{yzyz \ldots y(z)}_k + 
\underbrace{zyzy \ldots z(y)}_k 
\end{equation}
for any $k$.
Now let $f=(y+z)^{n-3}g$ for some word $g \neq 1$ beginning with $z$. 
There are two cases to consider:
$g=z$ and $g=zyh$ for some $h$ ($h=1$ is allowed). 
If $g=z$ equation \eqref{Binom} implies
(recall that $f$ is a cyclic word)
\begin{equation}\nonumber
\begin{cases}
f = 0 &\text{ if $n$ is odd } \\
2f = (y+z)^{n-2} = 0 &\text{ if $n$ is even }
\end{cases}
\end{equation}
If $g=zyh$ we have 
\begin{equation}\nonumber
f = \underbrace{(y)zyz \ldots yz}_{n-2} yh
\overset{\eqref{DnYZRelation}}{=}
- \underbrace{(z)yzy \ldots zy}_{n-2} yh
\overset{\eqref{DnYRelation}}{=} 0 \medspace ,
\end{equation}
which completes the proof for $D_n$ assuming $\charac k \neq 2$. 

\subsection{$E_6$-case}
We have 
\begin{gather}
y^3 = 0 \\ \label{E6ZRelation}
z^2 = 0 \\ \label{E6YZRelation}
(y+z)^{3} = y^2z + yzy + zy^2 + zyz = 0
\end{gather}
and $f = (y + z)^2 g = y^2 g$ for some word $g \neq 1$
beginning and ending with $z$. There are two cases:
$g=z$ and $g=zhz$ for some word $h$ beginning and ending with $y$.
If $g=z$ we have (recall that $f$ is a cyclic word)
\begin{equation}\nonumber
3f = 3 y^2 z = (y + z)^3 \overset{\eqref{E6YZRelation}}{=} 0 
\medspace .
\end{equation}
If $g = zhz$ we have
\begin{equation}\nonumber
f = y^2zhz \overset{\eqref{E6YZRelation}}{=} 
-yzyhz = - zyzyh 
\overset{\eqref{E6YZRelation}}{=} 
(y^2z + yzy + zy^2)yh = 0
\end{equation}
since $h$ begins and ends with $y$ and $y^3 =0$.
This completes the proof for $E_6$ assuming 
$\charac k \neq 2$, $3$ (one needs the assumption $\charac k \neq 2$
because $E_6$ contains subgraphs of type $D_n$).

\subsection{$E_7$-case}
We have 
\begin{gather} \label{E7YRelation}
y^3 = 0 \\ \label{E7ZRelation}
z^2 = 0 \\ \label{E7YZRelation}
(y+z)^{4} = y^2zy + yzy^2 + yzyz + zy^2z + zyzy = 0
\end{gather}
and $f = (y + z)^3 g = yzy g$ for some word $g \neq 1$
beginning and ending with $z$. There are two cases:
$g=z$ and $g=zhz$ for some word $h$ beginning and ending with $y$.
If $g=z$ we have (recall that $f$ is a cyclic word)
\begin{equation}\nonumber
2f = 2 yzyz \overset{(\ref{E7YRelation}, \ref{E7ZRelation})}{=} 
(y + z)^4 \overset{\eqref{E7YZRelation}}{=} 0 
\medspace .
\end{equation}
If $g = zhz$ we have
\begin{equation}\nonumber
f = yzyzhz \overset{\eqref{E7YZRelation}}{=} 
-y^2zyhz = -zy^2zyh
\end{equation}
since $h$ begins and ends with $y$ and $y^3 =0$.
So $f$ is a cyclic word in syllabi $zy$ and $zy^2$ having 
length greater than 5 and containing $zy^2$. 
We are going to prove that any such word is equal to $0$. 

Using \eqref{E7YZRelation} we get the following relations:
\begin{align}\label{E7FFS}
zyzyzy^2 = 0 \qquad &\text{ [ $z\eqref{E7YZRelation}y^2$ ] ,} \\
\label{E7FSS}
zy^2zyzy^2 + zyzy^2zy^2 = 0  \qquad
&\text{ [ $z\eqref{E7YZRelation}zy^2$ ] ,} \\ 
\label{E7SSS}
zy^2zy^2zy^2 = 0 \qquad &\text{ [ $zy^2\eqref{E7YZRelation}y^2$ ] } 
\medspace .
\end{align}
If $f$ contains only one $zy^2$ and has length greater than 5 it
vanishes because of \eqref{E7FFS}. If there are more than
one $zy^2$ we can use \eqref{E7FSS} to move syllabi of the type $zy$
together and then conclude from \eqref{E7FFS} that the cyclic 
word $f$ vanishes if it contains more than one such syllabi.
If $f$ contains one $zy$ then \eqref{E7FSS} implies that $2f =0$.
The relation \eqref{E7SSS} proves that any word having more than
two syllabi $zy^2$ and having no syllabi $zy$ vanishes, so the
only remaining case is $f = zy^2zy^2$, which is equal to $0$
because of (\ref{E7YZRelation}, \ref{E7YRelation}).
This completes the proof for $E_7$ assuming 
$\charac k \neq 2, 3$.

\subsection{$E_8$-case}
We have 
\begin{gather} \label{E8YRelation}
y^3 = 0 \\ \label{E8ZRelation}
z^2 = 0 
\end{gather}
\begin{multline}\label{E8YZRelation}
(y+z)^{5} = \\ =
y^2zy^2 + y^2zyz + yzy^2z + yzyzy + zy^2zy + zyzy^2 + zyzyz = 0
\end{multline}
and $f = (y + z)^4 g = y^2zy g + yzy^2 g$ for some word $g \neq 1$
beginning and ending with $z$. There are two cases:
$g=z$ and $g=zhz$ for some word $h$ beginning and ending with $y$.
If $g=z$ we have (recall that $f$ is a cyclic word)
\begin{equation}\nonumber
5f = 10 zyzy^2 \overset{(\ref{E8YRelation}, \ref{E8ZRelation})}{=} 
2(y + z)^5 \overset{\eqref{E8YZRelation}}{=} 0 
\medspace .
\end{equation}
If $g = zhz$ we have
\begin{equation}\nonumber
f = y^2zyzhz + yzy^2zhz = zy^2 zyzh + zyzy^2zh \medspace .
\end{equation}
Hence $f$ is a sum of cyclic words in syllabi $zy$ and $zy^2$ having 
length greater than 6 and containing $zy^2$. 
We are going to prove that any such word is equal to $0$. 

Using \eqref{E8YZRelation} we get the following relations:
\begin{align}\label{E8FFS}
zy^2zyzy + zyzy^2zy + zyzyzy^2 = 0 
\qquad &\text{ [ $z\eqref{E8YZRelation}y$ ] ,} \\ 
\label{E8FSS}
zy^2zyzy^2 + zyzy^2zy^2 = 0 
\qquad &\text{ [ $z\eqref{E8YZRelation}y^2$ ] ,} \\
\label{E8SSS}
zy^2zy^2zy^2 + zyzyzyzy^2 = 0
\qquad &\text{ [ $z\eqref{E8YZRelation}zy^2$ ] ,} \\
\label{E8SFFS}
zy^2zyzyzy^2 = 0
\qquad &\text{ [ $zy^2\eqref{E8YZRelation}y^2$ ] .}
\end{align}
Consider a cyclic word $w$ in syllabi $zy$ and $zy^2$ containing $zy^2$ 
and having length greater than 6.
If $w$ has only one $zy^2$ then
\eqref{E8FFS} implies that $3w=0$. If there are more than
one $zy^2$'s we can first use \eqref{E8SSS} to remove strings of $zy$'s
longer than $zyzy$. Now if there are two or more $zy$'s we can move
them together using \eqref{E8FSS} and 
then conclude from \eqref{E8SFFS} that the cyclic word $w$ vanishes.
If $w$ contains one $zy$ then \eqref{E8FSS} implies that $2w =0$.
So the remaining case is a cyclic word having only syllabi
of the type $zy^2$ and, moreover, having three or more such syllabi
(since the length is greater than 6). In this case we can use
\eqref{E8SSS} to create a substring $zyzyzyzy^2$ inside the word $w$. 
Now if $w = zyzyzyzy^2$ then $3w = 0$ because of 
\eqref{E8FFS}. Otherwise $w$ contains the string
$zy^2zyzyzyzy^2$ and so it vanishes because of
(\ref{E8FFS}, \ref{E8SFFS}).
This completes the proof for $E_8$ assuming 
$\charac k \neq 2,3,5$.

\section{Quiver Varieties}\label{VarietiesSection}

In this section the graph $G$ is of $ADE$ type and
$k$ is an algebraically closed field of characteristic zero.

\subsection{Quivers.}

Since we plan to discuss, in particular, Poisson structures
which are antisymmetric by nature we need to choose
an orientation of $G$ (i.e. make it into a quiver $Q$)
and use the quiver-based algebras $\Pi^{\lambda}_{\Omega}$
and $L_{\Omega}^{\lambda}$ (cf. \ref{QuiverAlgebras}). 
The graph being Dynkin (and hence bipartite) 
it is easy to see that 
$L_{\Omega}^{\lambda} \simeq L^{\lambda}$ and 
$\Pi_{\Omega}^{\lambda} \simeq \Pi^{\lambda}$ as
graded algebras. 

\subsection{Varieties.}
Now let us turn to geometry.
Given $I$-graded vector spaces $D$ and $V$,
Nakajima \cite{Nakajima98} considers the affine algebraic
variety (in fact a linear space) $N_{D,V}$ of all triples
$(x, p, q)$, where $x$ 
is a graded representation of the path algebra $k\overline{Q}$ 
in $V$, and $p: D \rightarrow V$ and 
$q: V \rightarrow D$ are graded linear maps. 
The variety $N_{D,V}$ is Poisson with respect to the bivector 
$\sum_{i \in I} \Tr(\partial_{p_i} \wedge \partial_{q_i}) +
\sum_{a \in H} \Tr (\partial_{x(a)} \wedge \partial_{x(a^*)})$.

Let $GL(V)$ denote the group of graded automorphisms of $V$
and $\mathfrak{gl}(V)$ be its Lie algebra. 
Then $GL (V)$ acts naturally on $N_{D,V}$ preserving the Poisson
structure and with the moment map 
$\mu: N_{D,V} \rightarrow \mathfrak{gl} (V)^*$ given by
\begin{equation}\nonumber 
\mu (x,p,q) = 
\sum_{a \in H} \bigl( x(a) x(a^*) -
x(a^*) x(a) \bigr) - \sum_{i \in I} p_i q_i =
x ( \theta_{\Omega} ) - pq
\medspace ,
\end{equation}
where we consider $\lambda \in k^I$ as the diagonal matrix
$\sum_{i \in I} \lambda_i \id_{V_i} \in \mathfrak{gl} (V)$ and
identify $\mathfrak{gl} (V)$ and $\mathfrak{gl} (V)^*$
via the trace form.
Let $\mathfrak{N}_{D,V, \lambda}$ be the symplectic quotient
at the value of the moment map equal to $\lambda$:
\begin{multline}\nonumber
\mathfrak{N}_{D,V, \lambda} = \mu^{-1} (\lambda) 
/ \negthickspace / 
GL(V) 
= \\ =
\{ (x,p,q) \in N_{D,V} \ | \ 
x (\theta_{\Omega}-\lambda) = pq \}
/ \negthickspace / 
GL(V) 
\medspace .
\end{multline}
Here $/ \negthickspace /$ denotes the affine quotient
(spectrum of the ring of invariant global functions).
The associated reduced scheme $\mathfrak{N}^{\text{red}}_{D,V, \lambda}$ 
was introduced by
Nakajima \cite{Nakajima98} who also constructed a resolution of
singularities of $\mathfrak{N}^{\text{red}}_{D,V, \lambda}$.

On the other hand, following Lusztig \cite{Lusztig03}*{2.3},
one can consider the affine scheme $\mathfrak{L}_{D, \lambda}$ of
graded representations of the algebra $L^{\lambda}_{\Omega}$ in
the $I$-graded vector space $D$ (it is proven by
Lusztig \cite{Lusztig03}*{Lemma 2.2} that $L^{0}_{\Omega}$ is 
finitely generated and the proof works for arbitrary $\lambda \in k^I$). 
Note that $\mathfrak{L}_{D, \lambda}$ is the scheme of \emph{all}
representations not of isomorphism classes.

One has a natural map 
\begin{equation}\nonumber
\vartheta: \mathfrak{N}_{D,V,\lambda}
\rightarrow \mathfrak{L}_{D,\lambda} 
\end{equation}
given by 
\begin{equation}\nonumber
\vartheta ((x,p,q)) (f) = q x(f) p
\medspace ,  
\end{equation}
where $f \in k\overline{Q}$ is a path in $\overline{Q}$.
The moment map condition $x (\theta_{\Omega} - \lambda) = pq$ 
implies that
\begin{multline}\nonumber
\vartheta ((x, p, q)) (f) \vartheta ((x,p,q)) (f') 
= \\ =
\vartheta ((x, p, q)) (f \cdot (\theta_{\Omega} - \lambda) 
\cdot f')= 
\vartheta ((x, p, q)) (f \circ f') 
\medspace .
\end{multline}
Hence the map $\vartheta$ is well-defined.

Lusztig proved that $\vartheta$ is
injective and finite \cite{Lusztig00}*{Theorem 4.7}
and that $\mathfrak{L}_{D, \lambda}$
is the union of images of maps $\vartheta$ for various $V$
\cite{Lusztig00}*{Lemma 4.12}.
The following Theorem is a refinement of the Lusztig's result.

\begin{theorem}\label{ImmersionTheorem}
The map $\vartheta$ is an immersion of Nakajima's 
scheme $\mathfrak{N}_{D,V,\lambda}$
as a subscheme of the scheme $\mathfrak{L}_{D,\lambda}$
of representations of $L^\lambda_{\Omega}$ in $D$.
\end{theorem}
\begin{proof} 
It is shown by Lusztig 
\cite{Lusztig98}*{Section 1} and it also follows from
the result of Le Bruyn and Procesi 
\cite{LeBruynProcesi90}*{Section 3} using a trick
of Crawley-Boevey 
\cite{CrawleyBoevey01}*{Remarks after the Introduction}
that the algebra $k [\mathfrak{N}_{D,V,\lambda}]$
of global functions on 
$\mathfrak{N}_{D,V,\lambda}$ is generated by matrix elements
of the operators $q x(f) p : D \rightarrow D$ and
traces of the operators $x(f) : V \rightarrow V$. 
The Theorem means that elements of the first type
(which come from $\mathfrak{L}_{D,\lambda}$ via $\vartheta^{-1}$)
already generate $k [\mathfrak{N}_{D,V,\lambda}]$. 
Let us prove it. Consider $\Tr x(f)$ for some path 
$f \in k\overline{Q}$. We can use Theorem \ref{MainLemma}
(more precisely, its variant for $\Pi^{\lambda}_{\Omega}$, or
the original Theorem together with an isomorphism
$\Pi^{\lambda}_{\Omega} \simeq \Pi^{\lambda}$)
to write $f$ as a sum 
$f = f_{\theta} + f_{[,]} + f_0$, where $f_{\theta}$ belongs
to the two sided ideal generated by $\theta_{\Omega}$,
$f_{[,]} \in [k\overline{Q}, k\overline{Q}]$, 
and $f_0$ has length zero.
Now $\Tr x (f_0)$ is constant on $\mathfrak{N}_{D,V,\lambda}$
(it depends only on the graded dimension of $V$), 
$\Tr x (f_{[,]}) = 0$, and the moment
map condition $x (\theta_{\Omega}) = pq$ implies that
$\Tr x (f_{\theta}) = \Tr {q x (f') p}$ for
some $f' \in k\overline{Q}$.
\end{proof}

\noindent\textbf{Remark}. A similar statement in the affine case would say
that the fibers of the map $\vartheta$ are related to the 
symmetric powers of Kleinian
singularities. That statement is not made precise in the present 
paper (note in particular that
$\mathfrak{L}_{D,\lambda}$ is infinite-dimensional
in the affine case).
Also Lusztig's algebra $L^{\lambda}_{\Omega}$ in the affine case
can be described via McKay correspondence. 
Namely it is isomorphic as a vector space 
to $u ( k \langle x,y \rangle \# k\Gamma ) u$, where
$k \langle x,y \rangle$ is the ring of noncommutative polynomials in
two variables, $\Gamma$ is a finite subgroup of $SL (2, k)$, and
$u$ is the sum of orthogonal idempotents in $k\Gamma$,
but the multiplication in $L^{\lambda}_{\Omega}$ 
is given by the insertion
of the element $(xy - yx - \lambda)$ with $\lambda$ considered as an
element of the center of $k\Gamma$ 
(cf. \eqref{LusztigProduct}).

\subsection{Poisson structures}\label{PoissonSubSection}

Recall the following facts about Nakajima's varieties
$\mathfrak{N}^{\text{red}}_{D,V,\lambda}$
(cf. \citelist{\cite{Nakajima98} \cite{Lusztig03}}):
\begin{itemize}
\item
$\{ \mathfrak{N}^{\text{red}}_{D,V,\lambda} \}$
form an inductive system of Poisson varieties 
with respect to embeddings of the $I$-graded vector spaces $V$;
\item
$\mathfrak{L}^{\text{red}}_{D,\lambda} \subset 
\bigcup_{V} \image (\vartheta: \mathfrak{N}_{D,V,\lambda}
\rightarrow \mathfrak{L}_{D,\lambda})$;
\item
$\mathfrak{N}^{\text{red}}_{D,V,\lambda}$ contains a smooth
open (possibly empty) symplectic subset
$\mathfrak{N}^{s}_{D,V,\lambda}$ (the set of stable points)
such that the complement of $\mathfrak{N}^{s}_{D,V,\lambda}$
in $\mathfrak{N}^{\text{red}}_{D,V,\lambda}$ is the union of
images of the natural embedding maps 
$\mathfrak{N}^{\text{red}}_{D,V',\lambda} \rightarrow
\mathfrak{N}^{\text{red}}_{D,V,\lambda}$ for 
$V' \subset V$;
\item 
there are finitely many non-isomorphic $V$'s such that
$\mathfrak{N}^{s}_{D,V,\lambda}$ is non-empty.
\end{itemize}

Therefore Theorem \ref{ImmersionTheorem} implies the following.

\begin{corollary}
The variety $\mathfrak{L}^{\text{red}}_{D,\lambda}$ 
of representations of $L^{\lambda}_{\Omega}$ in $D$ is
Poisson. Moreover it is a union of finitely
many symplectic leaves $\vartheta (\mathfrak{N}^s_{D,V,\lambda})$;
in particular, each symplectic leaf is the smooth locus 
of its closure.
\end{corollary}

\noindent\textbf{Remark.} Instead of using the inductive system of
varieties $\mathfrak{N}^{\text{red}}_{D,V,\lambda}$ to
define the Poisson structure on the variety 
$\mathfrak{L}^{\text{red}}_{D,\lambda}$ of representations
of the algebra $L^{\lambda}_{\Omega}$ one can, in the spirit of
non-commutative symplectic geometry 
(cf. \cites{Kontsevich93, CrawleyBoeveyEtingofGinzburg05}), 
define a (non-commutative)
Poisson bracket directly on $L^{\lambda}_{\Omega}$ as follows:
\begin{equation}\nonumber
\{ f' , f'' \} = \sum_{a \in H} 
\bigl( P_a (f',f'') - P_a (f'',f') \bigr) +
f' \cdot f'' - f'' \cdot f' \medspace ,
\end{equation}
where the path $P_a (f', f'')$ is obtained from the paths 
$f'$ and $f''$ by locating an edge $a$ in $f'$ and an edge
$a^*$ in $f''$ (if there are any), removing these edges and
gluing the resultant paths together and $\cdot$ is, as before,
the concatenation of paths. Then the Poisson structure
descends to the variety $\mathfrak{L}^{\text{red}}_{D, \lambda}$
of representations of $L^{\lambda}_{\Omega}$
in any graded vector space $D$. This universal 
construction produces the same Poisson bracket as the one
defined via the maps $\vartheta: \mathfrak{N}_{D,V,\lambda}
\rightarrow \mathfrak{L}_{D,\lambda}$.

\end{document}